\documentclass[12pt]{amsart}
\newcommand{\RP}{{\mathbb {RP}}}
\newcommand{\CP}{{\mathbb {CP}}}
\newcommand{\Z}{{\mathbb Z}}
\newcommand{\R}{{\mathbb R}}
\newcommand{\C}{{\mathbb C}}
\newcommand{\Q}{{\mathbb Q}}
\newcommand{\sm}{\setminus}

\newtheorem{theorem}{Theorem}
\newtheorem{example}{Example}
\newtheorem{proposition}{Proposition}

\newtheorem{remark}{Remark}

\newtheorem{lemma}{Lemma}

\begin{document}
\author{V.A.~Vassiliev}
\address{Steklov Mathematical Institute of Russian Academy of Sciences; National Research Institute --- Higher School of Economics} \email{vva@mi.ras.ru}
\thanks{Supported by Program ``Leading scientific schools'', grant No. NSh-5138.2014.1
and RFBR grant 13-01-00383}

\title{Homology of spaces of non-resultant polynomial systems in $\R^2$ and $\C^2$}
\date{\today}

\begin{abstract}
The {\em resultant} variety in the space of systems of homogeneous polynomials
of given degrees consists of such systems having non-trivial solutions. We
calculate the integer cohomology groups of all spaces of non-resultant systems
of polynomials $\R^2 \to \R$, and also the rational cohomology groups of
similar systems in $\C^2$.
\end{abstract}
\maketitle

\section{Introduction}
Given $n$ natural numbers $d_1 \ge d_2 \ge \dots \ge d_n$, consider the space
of all real homogeneous polynomial systems
\begin{equation}
\label{msyst}
\left\{
\begin{aligned}
a_{1,0} x^{d_1} + a_{1,1} x^{d_1-1}y + \dots + a_{1,{d_1}}y^{d_1}  \\
\dots \dots \dots \dots \dots \dots \dots \dots \dots \dots \dots  \\
a_{n,0} x^{d_n} + a_{n,1} x^{d_n-1}y + \dots + a_{n,{d_n}}y^{d_n} \\
\end{aligned}
\right.
\end{equation}
in two real variables $x, y$.

We will refer to this space as $\R^D$, $D = \sum_1^n (d_i+1)$. The {\em
resultant variety} $\Sigma \subset \R^D$ is the space of all systems having
non-zero solutions. $\Sigma$ is a  semialgebraic subvarety of codimension $n-1$
in $\R^D$. Below we calculate the cohomology group of its complement, $H^*(\R^D
\setminus \Sigma)$.

For the ``affine'' version of this problem (concerning the space of
non-resultant systems of polynomials $\R^1 \to \R^1$ with leading terms
$x^{d_i}$) see e.g. \cite{V}, \cite{Vfil}, \cite{Kozl}. A similar calculation
for spaces of real homogeneous polynomials in $\R^2$ without zeros of
multiplicity $\ge m$ was done in \cite{hompol}.

Also, we calculate below the rational cohomology groups of the complex analogs
$\C^D \setminus \Sigma_\C$ of all spaces $\R^D \setminus \Sigma$.

\section{Main results}

\subsection{Notation}
\label{nota} For any $p=1, 2, \dots, d_1$ denote by $N(p)$ the sum of all
numbers $d_i+1,$ $i=1, \dots, n,$ which are smaller than or equal to $p$, plus
$p$ times the number of those $d_i$ which are equal or greater than $p$. (In
other words, $N(p)$ is the area of the part of Young diagram $(d_1+1, \dots,
d_n+1)$ strictly to the left from the $(p+1)$th column.) Let the index
$\Upsilon(p)$ be equal to the number of even elements $d_i \ge p$ if $p$ is
even, and to the number of odd elements $d_i \ge p$ if $p$ is odd. By $\tilde
H^*(X)$ we denote the cohomology group reduced modulo a point.
${\overline{H}}_*(X)$ denotes the Borel--Moore homology group, i.e. the
homology group of the one-point compactification of $X$ reduced modulo the
added point.

\begin{theorem} \label{mthm}
If the space $\R^D \setminus \Sigma$ is non-empty $($i.e. either $n>1$ or $d_1$
is even$)$ then the group $\tilde H^*(\R^D \setminus \Sigma, \Z)$ is equal to
the direct sum of following groups:

\noindent A$)$ For any $p=1, \dots, d_3$,

if $\Upsilon(p)$ is even, then $\Z$ in dimension $N(p)-2p$ and $\Z$ in
dimension $N(p)-2p+1$,

if $\Upsilon(p)$ is odd, then only one group $\Z_2$ in dimension $N(p)-2p+1$;

\noindent B$)$ If $d_1-d_2$ is odd then only one additional summand $\Z$ in
dimension $D-d_1-d_2-2$. If $d_1-d_2$ is even, then additional summand
$\Z^{d_2-d_3+1}$ in dimension $D-d_1-d_2-1$ and $($if $d_2 \neq d_3)$ summand
$\Z^{d_2-d_3}$ in dimension $D-d_1-d_2-2$.
\end{theorem}

\begin{example} \label{mex} \rm Let $n=2$.
If $d_1$ and $d_2$ are of the same parity, then $\R^D \setminus \Sigma$
consists of $d_2+1$ connected components, each of which has the homology of a
circle. For the invariant separating systems from different components we can
take the index of the induced map of the unit circle $S^1 \subset \R^2$ into
$\R^2 \setminus 0$. This index can take all values of the same parity as $d_2$
from the segment $[-d_2,d_2]$. The 1-dimensional cohomology class inside any
component is just the rotation index of the image of a fixed point (say,
$(1,0)$) around the origin.

If $d_1$ and $d_2$ are of different parities, then the space $\R^D \setminus
\Sigma$ has the homology of a two-point set. The index separating its two
connected components can be calculated as the parity of the number of zeros of
the odd-degree polynomial of our non-resultant system, which lie in the
(well-defined) domain in $\RP^1$ where the even-degree polynomial is positive.
\end{example}

Now, let $\C^D$ be the space of all polynomial systems (\ref{msyst}) with {\em
complex} coefficients $a_{i,j}$, and $\Sigma_\C \subset \C^D$ the set of
systems having solutions in $\C^2 \setminus 0$.

\begin{theorem}
\label{mthmc} For any $n>1$ the group $H^*(\C^D \setminus \Sigma_{\C}, \Q)$ is
isomorphic to $\Q$ in dimensions $0, 2n-3, 2n-1$ and $4n-4$, and is trivial in
all other dimensions.
\end{theorem}

Consider also the space $\C^{d+1}$ of all complex homogeneous polynomials
$$ a_0 x^{d} + a_1 x^{d-1}y + \dots + a_dy^{d} $$
and $m$-{\em discriminant} $\Sigma_m$ in it, consisting of all polynomials
vanishing on some line with multiplicity $ \ge m$.

\begin{theorem}
\label{mthmc2} For any $m>1$ and $d \ge 2m$, the group $H^*(\C^{d+1} \setminus
\Sigma_m, \Q)$ is isomorphic to $\Q$ in dimensions $0, 2m-3, 2m-1$ and $4m-4$,
and is trivial in all other dimensions. For any $m>1$ and $d \in [m, 2m-1]$
this group is isomorphic to $\Q$ in dimensions $0, 2m-3, 2m-1$ and $2d-2$. and
is trivial in all other dimensions.
\end{theorem}

\section{Some preliminary facts}

Denote by $B(M,p)$ the {\em configuration space} of subsets of cardinality $p$
in the topological space $M$.

\begin{lemma} \label{lem1} For any natural $p,$ there is a locally trivial
fibre bundle $B(S^1,p) \to S^1$, whose fiber is homeomorphic to $\R^{p-1}$.
This fibre bundle is non-orientable if $p$ is even, and is orientable $($and
hence trivial$)$ if $p$ is odd. \quad $\Box$ \end{lemma}

Indeed, the projection of this fibre bundle can be realised as the product of
$p$ points of the unit circle in $\C^1$. The fibre of this bundle can be
identified in the terms of the universal covering $\R^p \to T^p$ with any
connected component of a hyperplane $\{x_1 + \dots + x_p = \mbox{const}\}$ from
which all affine planes given by $x_i = x_j + 2\pi k$, $i \neq j$, $k \in \Z$,
are removed. Such a component is convex and hence diffeomorphic to $\R^{p-1}$.
The assertion on orientability can be checked immediately. \hfill $\Box$
\medskip

Let us embed a manifold $M$ generically into the space $\R^T$ of a very large
dimension, and denote by $M^{*r}$ the union of all $(r-1)$-dimensional
simplices in $\R^T$, whose vertices lie in this embedded manifold (and the
``genericity'' of the embedding means that if two such simplices have a common
point in $\R^T$, then their minimal faces, containing this point, coincide).
\hfill $\Box$
\medskip

\begin{proposition}[C.~Caratheodory theorem, see also \cite{Vfil}] \label{carat}
For any $r \ge 1$, the space $(S^1)^{*r}$ is homeomorphic to $S^{2r-1}$. \quad
$\Box$ \end{proposition}

\begin{remark} \rm
This homeomorphism can be realized as follows. Consider the space $\R^{2r+1}$
of all real homogeneous polynomials $\R^2 \to \R^1$ of degree $2r$, the convex
cone in this space consisting of everywhere non-negative polynomials, and (also
convex) dual cone in the dual space $\widehat \R^{2r+1}$ consisting of linear
forms taking only positive values inside the previous cone. The intersection of
the boundary of this dual cone with the unit sphere in $\widehat \R^{2r+1}$ is
naturally homeomorphic to $(S^1)^{*r}$; on the other hand it is homeomorphic to
the boundary of a convex $2r$-dimensional domain.
\end{remark}

\begin{lemma}[see \cite{novikov}, Lemma 3]
\label{caratc} For any $r>1$ the group $H_*((S^2)^{*r}, \C)$ is trivial in all
positive dimensions. \hfill $\Box$ \end{lemma}

Consider the ``sign local system'' $\pm \Q$ over $B(\CP^1,p)$, i.e. the local
system of  groups with fiber $\Q$, such that the elements of
$\pi_1(B(\CP^1,p))$ defining odd (respectively, even) permutations of $p$
points in $\CP^1$ act in the fiber as multiplication by $-1$ (respectively,
$1$).

\begin{lemma}[see \cite{novikov}, Lemma 2]
\label{lem16} All groups $H_i(B(\CP^1,p);\pm \Q)$ with $p \ge 1$ are trivial
except only for $H_0(B(\CP^1,1), \pm \Q) \sim H_2(B(\CP^1,1), \pm \Q) \sim
H_2(B(\CP^1,2), \pm \Q) \sim \Q$. \hfill $\Box$
\end{lemma}

\section{Proof of Theorem \ref{mthm}}

Following \cite{A70}, we use the Alexander duality
\begin{equation}
\label{alex} \tilde H^i(\R^D \sm \Sigma) \simeq {\overline{H}}
_{D-i-1}(\Sigma),
\end{equation}
where ${\overline{H}} _*$ denotes the Borel---Moore homology.

\subsection{Simplicial resolution of $\Sigma$}
\label{sres}

To calculate the right-hand group in (\ref{alex}) we construct a {\it
resolution} of the space $\Sigma$. Let $\chi: \RP^1 \to \R^T$ be a generic
embedding, $T >> n$. For any system $\Phi=(f_1, \dots, f_n) \in \Sigma,$ not
equal identically to zero, consider the simplex $\Delta(\Phi)$ in $\R^T,$
spanned by the images $\chi(x_i)$ of all points $x_i \in \RP^1,$ corresponding
to all possible lines, on which the system $f$ has a common root. (The maximal
possible number of such lines is obviously equal to $d_1.$)

Further consider a subset in the direct product $\R^D \times \R^T$, namely the
union of all simplices of the form $\Phi \times \Delta(\Phi),$ $\Phi \in \Sigma
\sm 0$. This union is not closed: the set of its limit points, not belonging to
it, is the product of the point $0 \in \R^D$ (corresponding to the zero system)
and the union of all simplices in $\R^T,$ spanned by the images of no more than
$d_1$ different points of the line $\RP^1.$ By the Caratheodory theorem, the
latter union is homeomorphic to the sphere $S^{2d_1-1}.$ We can assume that our
embedding $\chi: \RP^1 \to \R^T$ is algebraic, and hence this sphere is
semialgebraic. Take a generic $2d_1$-dimensional semialgebraic disc in $\R^T$
with boundary at this sphere (e.g., the union of segments connecting the points
of this sphere with a generic point in $\R^T$) and add the product of the point
$0 \in \R^D$ and this disc to the previous union of simplices in $\R^D \times
\R^T$. The resulting set will be denoted by $\sigma$ and called a {\em
simplicial resolution} of $\Sigma$. \medskip

\begin{lemma} \label{lem2}
The obvious projection $\sigma \to \Sigma$ $($induced by the projection of
$\R^D \times \R^T$ onto the first factor$)$ is proper, and the induced map of
one-point compactifications of these spaces is a homotopy equivalence.
\end{lemma}

This follows easily from the fact that this projection is a stratified map of
semialgebraic spaces, and the preimage of any point $\bar \Sigma$ is
contractible, cf. \cite{V}, \cite{Vfil}. \quad $\Box$ \medskip

So, we can (and will) calculate the group ${\overline{H}} _*(\sigma)$ instead
of ${\overline{H}} _*(\Sigma)$.

\begin{remark} \rm
There is a more canonical construction of a simplicial resolution of $\Sigma$
in the terms of ``Hilbert schemes''. Namely, let $I_p$ be the space of all
ideals of codimension $p$ in the space of smooth functions $\RP^1 \to \R^1$,
supplied with the natural ``Grassmannian'' topology. It is easy to see that
$I_p$ is homeomorphic to the $p$th symmetric power $S^p(\RP^1) =
(\RP^1)^p/S(p),$ in particular it contains the configuration space $B(\RP^1,p)$
as an open dense subset. Consider the disjoint union of these $d_1$ spaces
$I_1, \dots, I_{d_1}$, augmented with the one-point set $I_0$ symbolizing the
zero ideal. The incidence of ideals makes this union a partially ordered set.
Consider the continuous order complex $\Xi_{d_1}$ of this poset, i.e. the
subset in the join $I_1 * \dots * I_{d_1} * I_0$ consisting of those simplices,
all whose vertices are incident to one another. For any polynomial system
$\Phi=(f_1, \dots, f_n) \in \R^D$ denote by $\Xi(\Phi)$ the subcomplex in
$\Xi_{d_1}$ consisting of all simplices, all whose vertices correspond to
ideals containing all polynomials $f_1, \dots, f_n$. The simplicial resolution
$\tilde \sigma \subset \Sigma \times \Xi_{d_1}$ is defined as the union of
simplices $\Phi \times \Xi(\Phi)$ over all $\Phi \in \Sigma$.

This construction is homotopy equivalent to the previous one. In particular,
the Caratheodory theorem has the following version: the continuous order
complex of the poset of all ideals of codimension $\le r$ in the space of
functions $S^1 \to \R^1$ is {\em homotopy equivalent} to $S^{2r-1}$.

However, this construction is less convenient for our practical calculations.
\end{remark}

The space $\sigma$ has a natural increasing filtration $F_1 \subset \dots
\subset F_{d_1+1} \equiv \sigma$: its term $F_p,$ $p \le d_1,$ is the closure
of the union of all simplices of the form  $\Phi \times \Delta(\Phi)$ over all
polynomial systems $\Phi$ having no more than $p$ lines of common zeros.
\medskip

\begin{lemma} \label{lem3} For any $p =1, \ldots, d_1,$ the
term $F_p \sm F_{p-1}$ of our filtration is the space of a locally trivial
fiber bundle over the configuration space $B(\RP^1,p),$ with fibers equal to
the direct product of an $(p-1)$-dimensional open simplex and an
$(D-N(p))$-dimensional real space. The corresponding bundle of open simplices
is orientable if and only if $p$ is odd $($i.e. exactly when the base
configuration space is orientable$)$, and the bundle of $(D-N(p))$-dimensional
spaces is orientable if and only if the index $\Upsilon(p)$ is even.

The last term $F_{d_1+1} \sm F_{d_1}$ of this filtration is homeomorphic to an
open $2d_1$-dimensional disc.
\end{lemma}

Indeed, to any configuration $(x_1, \ldots, x_p) \in B(\RP^1,p),$ $p \le d_1$,
there corresponds the direct product of the interior part of the simplex in
$\R^T,$ spanned by the images $\chi(x_i)$ of points of this configuration, and
the subspace in $\R^D,$ consisting of polynomial systems, having solutions on
corresponding $p$ lines in $\R^2.$ The codimension of the latter subspace is
equal exactly to $N(p)$. The assertion concerning the orientations can be
checked elementary. The description of  $F_{d_1+1} \sm F_{d_1}$ follows
immediately from the construction. \quad $\Box$
\medskip

Consider the spectral sequence $E_{p,q}^r,$ calculating the group
${\overline{H}} _*(\Sigma)$ and generated by this filtration. Its term
$E_{p,q}^1$ is canonically isomorphic to the group ${\overline{H}} _{p+q}(F_p
\sm F_{p-1}).$ By Lemma \ref{lem3}, its column $E_{p,*}^1,$ $p \le d_1,$ is as
follows. If $\Upsilon(p)$ is even, then it contains exactly two non-trivial
terms $E_{p,q}^1$, both isomorphic to $\Z$, for $q$ equal to $D-N(p)+p-1$ and
$D-N(p)+p-2$. If $\Upsilon(p)$ is odd, then it contains only one non-trivial
term $E_{p,q}^1$, isomorphic to $\Z_2$, for $q=D-N(p)+p-2$. Finally, the column
$E^1_{d_1+1,*}$ contains only one non-trivial element $E^1_{d_1+1,d_1-1} \sim
\Z$. \medskip

Before calculating the differentials and further terms $E^r$, $r >1$, let us
consider several basic examples.

\subsection{The case $n=1$}

\begin{figure}
\unitlength=0.80mm \special{em:linewidth 0.4pt} \linethickness{0.4pt} \mbox{
\begin{picture}(73.00,45.00)
\thicklines \put(10.00,10.00){\vector(1,0){63.00}} \put(10.00,10.00){\vector(0,1){35.00}} \thinlines
\put(13.00,6.00){\makebox(0,0)[cc]{\small $1$}} \put(19.00,6.00){\makebox(0,0)[cc]{\small $2$}}
\put(25.00,6.00){\makebox(0,0)[cc]{\small $3$}} \put(61.00,6.00){\makebox(0,0)[cc]{\small $d_1$}}
\put(70.00,6.00){\makebox(0,0)[cc]{$p$}} \put(5.00,33.00){\makebox(0,0)[cc]{\small $d_1$}}
\put(13.00,27.00){\makebox(0,0)[cc]{${\mathbb Z}$}} \put(13.00,33.00){\makebox(0,0)[cc]{${\mathbb Z}$}}
\put(19.00,27.00){\makebox(0,0)[cc]{${\mathbb Z}_2$}} \put(19.00,33.00){\makebox(0,0)[cc]{$0$}}
\put(25.00,27.00){\makebox(0,0)[cc]{${\mathbb Z}$}} \put(25.00,33.00){\makebox(0,0)[cc]{${\mathbb Z}$}}
\put(34.00,33.00){\makebox(0,0)[cc]{$\dots$}} \put(34.00,27.00){\makebox(0,0)[cc]{$\dots$}}
\put(43.00,27.00){\makebox(0,0)[cc]{${\mathbb Z}$}} \put(43.00,33.00){\makebox(0,0)[cc]{${\mathbb Z}$}}
\put(49.00,27.00){\makebox(0,0)[cc]{${\mathbb Z}_2$}} \put(49.00,33.00){\makebox(0,0)[cc]{$0$}}
\put(55.00,27.00){\makebox(0,0)[cc]{${\mathbb Z}$}} \put(55.00,33.00){\makebox(0,0)[cc]{${\mathbb Z}$}}
\put(61.00,27.00){\makebox(0,0)[cc]{${\mathbb Z}_2$}} \put(61.00,33.00){\makebox(0,0)[cc]{$0$}}
\put(67.00,27.00){\makebox(0,0)[cc]{${\mathbb Z}$}} \put(6.00,43.00){\makebox(0,0)[cc]{$q$}}
\end{picture}
}
\mbox{
\begin{picture}(73.00,45.00)
\thicklines \put(10.00,10.00){\vector(1,0){63.00}} \put(10.00,10.00){\vector(0,1){35.00}} \thinlines
\put(13.00,6.00){\makebox(0,0)[cc]{\small $1$}} \put(19.00,6.00){\makebox(0,0)[cc]{\small $2$}}
\put(25.00,6.00){\makebox(0,0)[cc]{\small $3$}} \put(61.00,6.00){\makebox(0,0)[cc]{\small $d_1$}}
\put(70.00,6.00){\makebox(0,0)[cc]{$p$}} \put(5.00,33.00){\makebox(0,0)[cc]{\small $d_1$}}
\put(13.00,33.00){\makebox(0,0)[cc]{$0$}} \put(13.00,27.00){\makebox(0,0)[cc]{${\mathbb Z}_2$}}
\put(19.00,27.00){\makebox(0,0)[cc]{${\mathbb Z}$}} \put(19.00,33.00){\makebox(0,0)[cc]{${\mathbb Z}$}}
\put(25.00,27.00){\makebox(0,0)[cc]{${\mathbb Z}_2$}} \put(25.00,33.00){\makebox(0,0)[cc]{$0$}}
\put(34.00,33.00){\makebox(0,0)[cc]{$\dots$}} \put(34.00,27.00){\makebox(0,0)[cc]{$\dots$}}
\put(43.00,27.00){\makebox(0,0)[cc]{${\mathbb Z}$}} \put(43.00,33.00){\makebox(0,0)[cc]{${\mathbb Z}$}}
\put(49.00,27.00){\makebox(0,0)[cc]{${\mathbb Z}_2$}} \put(49.00,33.00){\makebox(0,0)[cc]{$0$}}
\put(55.00,27.00){\makebox(0,0)[cc]{${\mathbb Z}$}} \put(55.00,33.00){\makebox(0,0)[cc]{${\mathbb Z}$}}
\put(61.00,27.00){\makebox(0,0)[cc]{${\mathbb Z}_2$}} \put(61.00,33.00){\makebox(0,0)[cc]{$0$}}
\put(67.00,27.00){\makebox(0,0)[cc]{${\mathbb Z}$}} \put(6.00,43.00){\makebox(0,0)[cc]{$q$}}
\end{picture}
} \caption{$E^1$ for $n=1$, $d_1$ even \hspace{1.2cm} and $n=1$, $d_1$ odd}
\label{fig1}
\end{figure}
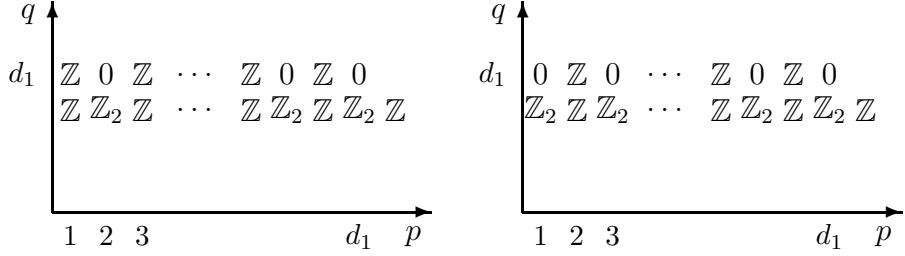

If our system consists of only one polynomial of degree $d_1$, then the term
$E^1$ of our spectral sequence looks as in Fig.~\ref{fig1}, in particular all
non-trivial groups $E_{p,q}^1$ lie in two rows $q=d_1$ and $d_1-1$.

\begin{lemma} \label{lem4}
If $n=1$ then in both cases of even or odd $d_1$, all possible horizontal
differentials $\partial_1: E_{p,d_1-1}^1 \to E_{p-1,d_1-1}^1$ of the form $\Z
\to \Z_2$, $p=d_1+1, d_1-1, d_1-3, \dots$ are epimorphic, and all differentials
$\partial_2: E_{p,d_1-1}^2 \to E_{p-2,d_1}^2$ of the form $\Z \to \Z$,
$p=d_1+1, d_1-1, d_1-3, \dots$ are isomorphisms. In particular, the unique
surviving term $E_{p,q}^3$ for the ``even'' spectral sequence is $E_{1,d_1-1}^3
\sim \Z$, and for the ``odd'' one it is $E_{2,d_1-1}^3 \sim \Z$.
\end{lemma}

Indeed, in both cases we know the answer. In the ``odd'' case the discriminant
coincides with entire $\R^D = \R^{d_1+1}$. In the ``even'' one its complement
consists of two contractible components, so that ${\overline{H}} _*(\Sigma) =
\Z$ in dimension $d_1$ and is trivial in all other dimensions. Therefore all
terms $E_{p,q}$ with $p+q$ not equal to $d_1+1$ (respectively, $d_1$) in the
odd (even) case should die on some step. \hfill $\Box$ \medskip

\subsection{The case $n=2$}

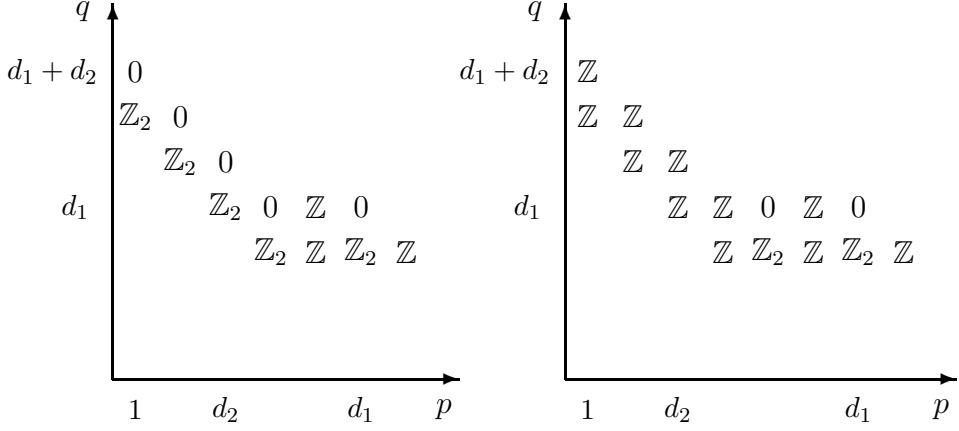
\begin{figure}
\unitlength=1.00mm
\special{em:linewidth 0.4pt}
\linethickness{0.4pt}
\mbox{
\begin{picture}(56.00,60.00)
\thicklines \put(10.00,10.00){\vector(1,0){46.00}} \put(10.00,10.00){\vector(0,1){50.00}} \thinlines
\put(13.00,6.00){\makebox(0,0)[cc]{\small $1$}} \put(25.00,6.00){\makebox(0,0)[cc]{\small $d_2$}}
\put(43.00,6.00){\makebox(0,0)[cc]{\small $d_1$}} \put(54.00,6.00){\makebox(0,0)[cc]{$p$}}
\put(5.00,33.00){\makebox(0,0)[cc]{\small $d_1$}} \put(8.00,51.00){\makebox(0,0)[rc]{\small $d_1+d_2$}}
\put(13.00,51.00){\makebox(0,0)[cc]{$0$}} \put(13.00,45.00){\makebox(0,0)[cc]{${\mathbb Z}_2$}}
\put(19.00,45.00){\makebox(0,0)[cc]{$0$}} \put(19.00,39.00){\makebox(0,0)[cc]{${\mathbb Z}_2$}}
\put(25.00,39.00){\makebox(0,0)[cc]{$0$}} \put(25.00,33.00){\makebox(0,0)[cc]{${\mathbb Z}_2$}}
\put(31.00,33.00){\makebox(0,0)[cc]{$0$}} \put(31.00,27.00){\makebox(0,0)[cc]{${\mathbb Z}_2$}}
\put(37.00,27.00){\makebox(0,0)[cc]{${\mathbb Z}$}} \put(37.00,33.00){\makebox(0,0)[cc]{${\mathbb Z}$}}
\put(43.00,33.00){\makebox(0,0)[cc]{$0$}} \put(43.00,27.00){\makebox(0,0)[cc]{${\mathbb Z}_2$}}
\put(49.00,27.00){\makebox(0,0)[cc]{${\mathbb Z}$}} \put(6.00,59.00){\makebox(0,0)[cc]{$q$}}
\end{picture}
}
\mbox{
\begin{picture}(62.00,60.00)
\thicklines \put(10.00,10.00){\vector(1,0){52.00}} \put(10.00,10.00){\vector(0,1){50.00}} \thinlines
\put(13.00,6.00){\makebox(0,0)[cc]{\small $1$}} \put(25.00,6.00){\makebox(0,0)[cc]{\small $d_2$}}
\put(49.00,6.00){\makebox(0,0)[cc]{\small $d_1$}} \put(60.00,6.00){\makebox(0,0)[cc]{$p$}}
\put(5.00,33.00){\makebox(0,0)[cc]{\small $d_1$}} \put(8.00,51.00){\makebox(0,0)[rc]{\small $d_1+d_2$}}
\put(13.00,51.00){\makebox(0,0)[cc]{${\mathbb Z}$}} \put(13.00,45.00){\makebox(0,0)[cc]{${\mathbb Z}$}}
\put(19.00,45.00){\makebox(0,0)[cc]{${\mathbb Z}$}} \put(19.00,39.00){\makebox(0,0)[cc]{${\mathbb Z}$}}
\put(25.00,39.00){\makebox(0,0)[cc]{${\mathbb Z}$}} \put(25.00,33.00){\makebox(0,0)[cc]{${\mathbb Z}$}}
\put(31.00,33.00){\makebox(0,0)[cc]{${\mathbb Z}$}} \put(31.00,27.00){\makebox(0,0)[cc]{${\mathbb Z}$}}
\put(37.00,33.00){\makebox(0,0)[cc]{$0$}} \put(37.00,27.00){\makebox(0,0)[cc]{${\mathbb Z}_2$}}
\put(43.00,33.00){\makebox(0,0)[cc]{${\mathbb Z}$}} \put(43.00,27.00){\makebox(0,0)[cc]{${\mathbb Z}$}}
\put(49.00,33.00){\makebox(0,0)[cc]{$0$}} \put(49.00,27.00){\makebox(0,0)[cc]{${\mathbb Z}_2$}}
\put(55.00,27.00){\makebox(0,0)[cc]{${\mathbb Z}$}} \put(6.00,59.00){\makebox(0,0)[cc]{$q$}}
\end{picture}
} \caption{$E^1$ for $n=2$, $d_1-d_2$ odd \hspace{1.2cm} and $n=2$, $d_1-d_2$
even} \label{fig2}
\end{figure}

There are two very different situations depending on the parity of $d_1-d_2$.
In Fig.~\ref{fig2} we demonstrate these situations in two particular cases,
$(d_1,d_2)=(6,3)$ and $(7,3)$. However, the general situation is essentially
the same, namely the following is true.

If $n=2$ and $d_1-d_2$ is odd, then all indices $\Upsilon(p)$, $p=1, \dots,
d_2+1$, are odd, and hence all non-trivial groups $E_{p,q}^1$ with such $p$ lie
on the line $\{p+q=d_1+d_2\}$ only and are equal to $\Z_2$.

If $n=2$ and $d_1-d_2$ is even, then all indices $\Upsilon(p)$, $p=1, \dots,
d_2+1$, are even, and hence all non-trivial groups $E_{p,q}^1$ with such $p$
lie on two lines $\{p+q=d_1+d_2\}$, $\{p+q=d_1+d_2+1\}$, and are all equal to
$\Z$.

In both cases, all groups $E_{p,q}^1$ with $p > d_2$ are the same as in the
case $n=1$ with the same $d_1$. Moreover, the differentials $\partial_1$ and
$\partial_2$ between these groups also are the same as for $n=1$, therefore all
of these groups die at the instant $E^3$ except for $E_{d_2+1,d_1-1}^3 \sim \Z$
for even $d_1-d_2$, and $E_{d_2+2,d_1-1}^3 \sim \Z$ for odd $d_1-d_2$.

In the case of even $d_1-d_2$ all other differentials between the groups
$E_{p,q}^r$ are trivial, because otherwise the group ${\overline{H}} ^0(\R^d
\setminus \Sigma)$ would be smaller than $\Z^{d_2}$, in contradiction to
$d_2+1$ different components of this space indicated in Example \ref{mex}.

On contrary, if $d_1-d_2$ is odd, then all differentials $d_r:E_{d_2+2,d_1-1}^r
\to E_{d_2+2-r, d_1-2+r}^r$, $r=1, \dots, d_1-d_2+1$ are epimorphic just
because the integer cohomology group of the topological space $\R^D \setminus
\Sigma$ cannot have non-trivial torsion subgroup in dimension 1. Therefore the
unique nontrivial group $E_{p,q}^\infty$ in this case is
$E_{d_2+2,d_1-1}^\infty \sim \Z$.

These considerations prove Main Theorem in the case $n=2$.

\subsection{The general case}

Now suppose that our systems (\ref{msyst}) consist of $n\ge 3$ polynomials. Let
again $\sigma$ be the simplicial resolution of the corresponding resultant
variety, constructed in \S \ref{sres}, and $\sigma'$ be the simplicial
resolution of the resultant variety for $n=2$ and the same $d_1$ and $d_2$. The
parts $\sigma \setminus F_{d_3}(\sigma)$ and $\sigma' \setminus
F_{d_3}(\sigma')$ of these resolutions are canonically homeomorphic to one
another as filtered spaces. In particular, $E_{p,q}^1(\sigma) =
E_{p,q}^1(\sigma')$ if $p
> d_3$, and $E_{p,q}^r(\sigma) = E_{p,q}^r$ if $p \ge d_3+r$. All non-trivial
terms $E_{p,q}^r(\sigma)$ with $p \le d_3$ are placed in such a way, that no
non-trivial differentials $\partial_r$ can act between these terms, as well as
no differentials can act to these terms from the cells $E_{p,q}^r$ with $p >
d_3$, which have survived the differentials between these cells, described in
the previous subsection.

Therefore the final term $E_{p,q}^\infty (\sigma)$ coincides with
$E_{p,q}^1(\sigma)$ in the domain $\{p \le d_3\}$, and coincides with the term
$E_{p,q}^\infty(\sigma')$ of the truncated spectral sequence calculating the
Borel-Moore homology of $\sigma' \setminus F_{d_3}(\sigma')$ in the domain $\{p
> d_3\}$. This terminates the proof of Theorem \ref{mthm}. \hfill $\Box$

\section{Proof of Theorems \ref{mthmc}, \ref{mthmc2}}

The simplicial resolution $\sigma_{\C}$ of $\Sigma_{\C}$ appears in the same
way as its real analog $\sigma$ in the previous section. It also has a natural
filtration $F_1 \subset \dots \subset F_{d_1+1} \equiv \sigma_{\C}$. For $p \in
[1, d_1]$ its term $F_p \setminus F_{p-1}$ is fibered over the configuration
space $B(\CP^1,p)$; its fiber over a configuration $(x_1, \dots x_p)$ is equal
to the product of the space $\C^{D-N(p)}$ (consisting of all complex systems
(\ref{msyst}) vanishing at all lines corresponding to the points of this
configuration) and the $(p-1)$-dimensional simplex, whose vertices correspond
to the points of the configuration. In particular, our spectral sequence
calculating  rational Borel-Moore homology of $\sigma_{\C}$ has $E^1_{p,q} \sim
{\overline{H}} _{q-2(D-N(p))+1}(B(\CP^1,p); \pm \Q)$ for such $p$. By Lemma
\ref{lem16}, only the following such groups are non-trivial: $E^1_{1,2(D-n)-1}
\sim \Q$, $E^1_{1,2(D-n)+1} \sim \Q$, and (if $d_1 >1$) $E^1_{2, 2(D-2n)+1}$.

The last term $F_{d_1+1} \setminus F_{d_1}$ is homeomorphic to the cone over
the $d_1$-th self-join $(\CP^1)^{*d_1}$ with the base of this cone removed (as
it belongs to $F_{d_1}$). Therefore by Lemma \ref{caratc} the column
$E^1_{d_1+1,*}$ is trivial if $d_1>1$ and contains unique non-trivial group
$E^1_{2,1} \sim \Q$ if $d_1=1$.

So, in any case the first leaf $E^1$ of our spectral sequence has only three
non-trivial terms $E^1_{1,2(D-n)-1} \sim \Q$, $E^1_{1,2(D-n)+1} \sim \Q$, and
$E^1_{2, 2(D-2n)+1}$. The differentials in it are obviously trivial, therefore
the group ${\overline{H}} _*(\sigma)$ has three non-trivial terms in dimensions
$2(D-n)$, $2(D-n)+2$, and $2(D-2n)+3$. By Alexander duality in the space $\C^D$
this gives us three groups $\tilde H^{2n-3} \sim \Q$, $\tilde H^{2n-1} \sim \Q$
and $\tilde H^{4n-4}\sim \Q$ and zero in all other dimensions. \hfill $\Box$

Theorem \ref{mthmc2} can be proved in exactly the same way. \hfill $\Box$


\begin{thebibliography}{99}

\bibitem{A70} V.I.~Arnold, {\it On some topological invariants of algebraic functions,} Trans. Moscow Math. Soc. 1970, V.~21, p.~30--52.

\bibitem{A89} V.I.~Arnold, {\it Spaces of functions with
moderate singularities,}
Funct. Anal. and its Appl. 23:3, (1989), p.~1--10.

\bibitem{CCMM} F.R.~Cohen, R.L.~Cohen, B.M.~Mann, R.J.~Milgram, {\it The topology of the space of rational functions and divisors of surfaces,} Acta Math., {\bf 166}, 1991, 163--221.

\bibitem{Kozl} A.~Kozlowski and R.~Yamaguchi, {\it Topology of complements of discriminants and resultants,}
Journal of the Mathematical Society of Japan, 52:4 (2000), 949-959.

\bibitem{V} V.A.~Vassiliev, {\it Complements of discriminants of smooth maps: topology and applications. Revised ed.}, Transl. of Math. Monographs, v.~98, AMS, Providence RI., 1994.

\bibitem{Vfil} V.A.~Vassiliev, {\it Topology of complements of discriminants,} Moscow,  Phasis, 1997, 552 p. (in Russian).

\bibitem{hompol} V.A.~Vassiliev, {\it Homology of spaces of homogeneous polynomials in $\R^2$ without multiple zeros}, Proc. Steklov Mathematical Institute, {\bf 221}, 1998, 143--148.

\bibitem{novikov} V.A.~Vassiliev, {\it How to calculte the homology of spaces of nonsingular algebraic projective hypersurfaces} Proc. Steklov Mathematical Institute, {\bf 225}, 1999, 121--140.


\end{thebibliography}
\end{document}